\documentclass{amsart}
\usepackage{url}
\usepackage{latexsym}
\usepackage{lscape} 
\usepackage{pdflscape}
\usepackage{float}
\usepackage{amssymb}
\usepackage{tikz}

\newtheorem{thm}[equation]{Theorem}

\newtheorem{exes}[equation]{Examples}
\newtheorem{lem}[equation]{Lemma}
\newtheorem{cor}[equation]{Corollary}

\newtheorem{assump}[equation]{Assumption}
\theoremstyle{remark}
\newtheorem{rem}[equation]{Remark}

\theoremstyle{definition}

\numberwithin{equation}{section}

\newcommand{\End}{\text{End}}
 
\newcommand{\Hom}{\text{Hom}}

\newcommand{\Ad}{\text{Ad}}
\newcommand{\ad}{\text{ad}}
\newcommand{\htil}{\tilde{h}}
\newcommand{\etil}{\tilde{e}}
\newcommand{\ftil}{\tilde{f}}

\newcommand{\fg}{{\mathfrak g}}

\newcommand{\fh}{{\mathfrak h}}
\newcommand{\fk}{{\mathfrak k}}
\newcommand{\fl}{{\mathfrak l}}

\newcommand{\fo}{{\mathfrak o}}
\newcommand{\fq}{{\mathfrak q}}

\newcommand{\fs}{{\mathfrak s}}

\newcommand{\Cal}{\mathcal}

\newcommand{\be}{\begin{equation}}
\newcommand{\beu}{\begin{equation*}}

\newcommand{\IP}[2]{\langle#1\,, #2\rangle}     

\usepackage{graphicx,color,cancel}

\newcommand\eps{{\varepsilon}}
\newcommand\la{{\lambda}}

\newcommand\twedge{\textstyle{\bigwedge}}

\def\dim{\mathop{\hbox {dim}}\nolimits}
\def\Ad{\mathop{\hbox {Ad}}\nolimits}
\def\ad{\mathop{\hbox {ad}}\nolimits}

\def\deg{\mathop{\hbox {deg}}\nolimits}

\def\End{\mathop{\hbox {End}}\nolimits}
\def\Hom{\mathop{\hbox {Hom}}\nolimits}

\def\ker{\mathop{\hbox{ker}}\nolimits}

\def\deg{\mathop{\hbox{deg}}\nolimits}
\def\id{\mathop{\hbox{id}}\nolimits}

\def\dim{\mathop{\hbox {dim}}\nolimits}

\newcommand{\lan}{\langle}
\newcommand{\ran}{\rangle}
\newcommand{\lara}{\langle\,,\rangle}

\newcommand{\pf}{\begin{proof}}
\newcommand{\epf}{\end{proof}}
\newcommand{\eq}{\begin{equation}}
\newcommand{\eeq}{\end{equation}}
\newcommand{\eqn}{\begin{equation*}}
\newcommand{\eeqn}{\end{equation*}}
\newcommand\bed{\begin{definition}}
\newcommand\ebed{\end{definition}}
\newcommand\bethm{\begin{theorem}}
\newcommand\ebethm{\end{theorem}}
\newcommand\bealigned{\begin{aligned}}
\newcommand\ebealigned{\end{aligned}}

\newcommand{\frg}{\mathfrak{g}}
\newcommand{\frh}{\mathfrak{h}}
\newcommand{\frk}{\mathfrak{k}}
\newcommand{\frl}{\mathfrak{l}}

\newcommand{\frq}{\mathfrak{q}}

\newcommand{\frs}{\mathfrak{s}}

\newcommand{\frsl}{\mathfrak{sl}}

\newcommand{\frso}{\mathfrak{so}}

\newcommand{\bbC}{\mathbb{C}}

\newcommand{\bbR}{\mathbb{R}}

\newcommand{\bbZ}{\mathbb{Z}}

\newcommand{\even}{\operatorname{even}}
\newcommand{\odd}{\operatorname{odd}}

\begin{document}

\sloppy
\title[Representation theoretic embedding of twisted Dirac operators]
{Representation theoretic embedding\\ of twisted Dirac operators}

\author{S. Mehdi}
\author{P. Pand\v zi\' c}
\address{Institut Elie Cartan de Lorraine, CNRS - UMR 7502, Universit\'e de Lorraine, France}
\email{salah.mehdi@univ-lorraine.fr}
\address{Department of Mathematics, Faculty of Science, University of Zagreb, Croatia}
\email{pandzic@math.hr}
\thanks{P.~Pand\v zi\'c was supported by the QuantiXLie  Center of Excellence, a project 
cofinanced by the Croatian Government and European Union through the European Regional Development Fund - the Competitiveness and Cohesion Operational Programme 
(KK.01.1.1.01.0004).}

\begin{abstract}
{\it Let $G$ be a non-compact connected semisimple real Lie group with finite center. Suppose $L$ is a non-compact connected closed subgroup of $G$ acting transitively on a symmetric space $G/H$ such that $L\cap H$ is compact. We study the action on $L/L\cap H$ of a Dirac operator $D_{G/H}(E)$ acting on sections of an $E$-twist of the spin bundle over $G/H$. As a byproduct, in the case of $(G,H,L)=(SL(2,{\mathbb R})\times SL(2,{\mathbb R}),\Delta(SL(2,{\mathbb R})\times SL(2,{\mathbb R})),SL(2,{\mathbb R})\times SO(2))$, we identify certain representations of $L$ which lie in the kernel of $D_{G/H}(E)$.} 
\end{abstract}

\maketitle

\section{Introduction}\label{Intro}

In a recent paper \cite{MO20} and an ongoing project \cite{MO21}, Mehdi and Olbrich study representation theoretic features of the spectrum of locally symmetric spaces arising as follows. Let $G$ be a non-compact connected semisimple real Lie group with finite center and $H$ a non-compact connected closed subgroup of $G$ so that 
$G/H$ is a symmetric space with respect to an involution $\sigma$ of $G$. Suppose there exists a non-compact connected real closed subgroup $L$ of $G$ which acts transitively and cocompactly on $G/H$, i.e., $G/H\simeq L/L\cap H$ with $L\cap H$ compact. In particular, if $D(G/H)$ (resp. $D(L/L\cap H)$) denotes the algebra of 
$G$-invariant (resp. $L$-invariant) differential operators on $G/H$ (resp. $L/L\cap H$), one gets an embedding of algebras 
\begin{equation}\label{introembed}
\imath:D(G/H)\hookrightarrow D(L/L\cap H).
\end{equation}
Now pick a discrete subgroup $\Gamma$ in $G$ which is contained in $L$ and such that $\Gamma\backslash G/H$ is a smooth compact manifold, i.e., a compact Clifford-Klein form of $G/H$. Using (\ref{introembed}), Mehdi and Olbrich describe the joint spectral decomposition of the (commutative) algebra $D(G/H)$ on the Hilbert space $L^2(\Gamma\backslash G/H)$ of square integrable functions on $\Gamma\backslash G/H$ in terms of the spectrum of $D(L/L\cap H)$ on $L^2(\Gamma\backslash L/L\cap H)\simeq L^2(\Gamma\backslash L)^{L\cap H}$. In particular, they prove the $L$-admissibility of certain $G$-representations involved in the spectrum. 

In this paper, we extend the embedding (\ref{introembed}) to the case of differential operators acting on sections of spin bundles over $G/H$ and $L/L\cap H$ twisted by a finite-dimensional representation $E$. Then we show that the image, under this embedding, of the Dirac operator $D_{G/H}(E)$ on $G/H$ splits as combination of both geometric and algebraic Dirac operators attached to the homogeneous spaces $L/L\cap H$, $L/L\cap K$, $L\cap K/L\cap H$ and $H/H\cap K$ (here $K$ is a maximal compact subgroup of $G$ fixed by a Cartan involution $\theta$, and $H$ and $L$ are assumed to be $\theta$-stable). Note that while most of these spaces are symmetric, $L/L\cap H$ is not, and for this space we use the cubic Dirac operator defined in \cite{G} and \cite{K}.

Finally, we use this splitting formula to relate the kernels of the various Dirac operators and we identify certain representations of $L$ occuring in the kernel of $D_{G/H}(E)$ when $G=SL(2,{\mathbb R})\times SL(2,{\mathbb R})$, $H=\Delta(SL(2,{\mathbb R})\times SL(2,{\mathbb R}))$ is the diagonal in $G$ and $L=SL(2,{\mathbb R})\times SO(2)$.

The paper is organized as follows. In Section \ref{triples} we define the triples $(G,H,L)$ and we collect their main features. In Section \ref{Spin} we review some facts about Clifford algebras and spin modules. In Section \ref{Dirac}, we describe the geometric and the algebraic Dirac operators, as well as the various bundles and sections on which they act. In Section \ref{embedding}, we compute the transfer of the cubic Dirac operator from $G/H$ to $L/L\cap H$ in terms of ``smaller'' Dirac operators. Finally, in Section \ref{sec sl2}, we use this transfer formula when $(G,H,L)=(SL(2,{\mathbb R})\times SL(2,{\mathbb R}), \Delta(SL(2,{\mathbb R})\times SL(2,{\mathbb R})),SL(2,{\mathbb R})\times SO(2))$
to relate the kernels of various Dirac operators. Moreover, we identify some of the representations of $L$ involved in the kernel of the twisted cubic Dirac operator on $G/H$.

\section{Transitive triples}\label{triples}

Let $G$ be a non-compact connected semisimple real Lie group with finite center and Lie algebra ${\mathfrak g}$ and let $\lan,\ran $ be the Killing form on ${\mathfrak g}$. Fix a Cartan involution $\theta$ on $G$ and let $K$ be the corresponding maximal compact subgroup of $G$ with Lie algebra ${\mathfrak k}$. The associated Cartan decomposition of ${\mathfrak g}$ is 
\begin{equation*}
{\mathfrak g}={\mathfrak k}\oplus{\mathfrak s},\quad\text{ with }[{\mathfrak k},{\mathfrak s}]\subset{\mathfrak s}\text{ and }[{\mathfrak s},{\mathfrak s}]\subset{\mathfrak k}.
\end{equation*}
Let $H$ be a non-compact connected semisimple closed subgroup of $G$ with Lie algebra ${\mathfrak h}$ such that the homogeneous space $G/H$ is a symmetric space with respect to some involution $\sigma$. There is therefore a decomposition of ${\mathfrak g}$: 
\begin{equation*}
{\mathfrak g}={\mathfrak h}\oplus{\mathfrak q},\quad\text{ with }[{\mathfrak h},{\mathfrak q}]\subset {\mathfrak q}\text{ and }[{\mathfrak q},{\mathfrak q}]\subset {\mathfrak h}.
\end{equation*}
Recall the following invariance properties of the Killing form. For all $X$, $Y$ and $Z$ in ${\mathfrak g}$, one has:
\begin{itemize}
\item[] $\lan\theta(X),\theta(Y)\ran =\lan X,Y\ran $
\item[] $\lan\sigma(X),\sigma(Y)\ran =\lan X,Y\ran $
\item[] $\lan \ad(X)(Z),Y\ran =-\lan Z,\ad(X)Y\ran $
\end{itemize}
Moreover, the restrictions of $\lan ,\ran $ to ${\mathfrak k}\times {\mathfrak k}$, ${\mathfrak s}\times {\mathfrak s}$, ${\mathfrak h}\times {\mathfrak h}$ and ${\mathfrak q}\times {\mathfrak q}$ are non-degenerate so that
\begin{equation*}
{\mathfrak k}\perp{\mathfrak s}\text{ and }{\mathfrak h}\perp{\mathfrak q}\text{ with respect to }\lan ,\ran .
\end{equation*}
We assume that $\theta$ and $\sigma$ commute, so that ${\mathfrak k}$ and ${\mathfrak s}$ are $\sigma$-stable:
\begin{equation*}
{\mathfrak k}=({\mathfrak k}\cap{\mathfrak h})\oplus({\mathfrak k}\cap{\mathfrak q})\text{ and }{\mathfrak s}=({\mathfrak s}\cap{\mathfrak h})\oplus({\mathfrak s}\cap{\mathfrak q})
\end{equation*}

Next, let $L$ be a non-compact connected semisimple closed subgroup of $G$ with Lie algebra ${\mathfrak l}$ such that 
\begin{itemize}\label{transtriples}
\item[(i)] $L$ is reductively embedded in $G$,
\item[(ii)] $L$ acts transitively on $G/H$,
\item[(iii)] $L\cap H$ is compact.
\end{itemize}
We will refer to triples $(G,H,L)$ satisfying (i),(ii) and (iii) as {\it transitive triples}. $L$ need not be $\sigma$-stable in general, however we may assume that $L$ is 
$\theta$-stable:
\begin{equation*}
{\mathfrak l}=({\mathfrak l}\cap{\mathfrak k})\oplus({\mathfrak l}\cap{\mathfrak s}).
\end{equation*}
The restriction to ${\mathfrak l}$ of the Killing form $\lan ,\ran $ of ${\mathfrak g}$ remains non-degenerate and 
\begin{equation*}
{\mathfrak g}={\mathfrak h}+{\mathfrak l}\text{ and }{\mathfrak q}\cap{\mathfrak l}^{\perp}=\{0\}
\end{equation*}
where ${\mathfrak l}^\perp$ is the orthogonal of ${\mathfrak l}$ in ${\mathfrak g}$ with respect to the Killing form. Transitive triples were classified in the late 1960s 
by Oni\v s\v cik in a more general setting \cite{Oni69}. A complete list can be found in \cite{KY05} or \cite{MO21}.
\begin{exes}\label{listtriples}\text{}\\
\begin{itemize}
\item[(1a)] $G=G^{\prime}\times G^{\prime}$, $H=\Delta (G^{\prime}\times  G^{\prime})$, $L=G^{\prime}\times\{e\}$.
\item[(1b)] $G=G^{\prime}\times G^{\prime}$, $H=\Delta (G^{\prime}\times G^{\prime})$, $L=G^{\prime}\times K^{\prime}$\\
where $G^\prime$ is a non-compact connected semisimple real Lie group with finite center and $K^\prime$ is a maximal compact subgroup of $G^{\prime}$.
\item[(2)] $G=SO_{e}(2,2n)$, $H=SO_{e}(1,2n)$, $L=U(1,n)$, $n\geq 1$.
\item[(3)] $G=SO_{e}(2,2n)$, $H=U(1,n)$, $L=SO_{e}(1,2n)$, $n\geq 1$.
\item[(4)] $G=SO_{e}(4,4n)$, $H=SO_{e}(3,4n)$, $L=Sp(1,n)$, $n\geq 1$.
\item[(5)] $G=SU(2,2n)$, $H=SU(1,2n)$, $L=Sp(1,n)$, $n\geq 1$.
\item[(6)] $G=SU(2,2n)$, $H=Sp(1,n)$, $L=SU(1,2n)$, $n\geq 1$.
\item[(7)] $G=SO_{e}(8,8)$, $H=SO_{e}(7,8)$, $L=Spin_{e}(1,8)$.
\item[(8)] $G=SO_{e}(4,4)$, $H=SO_{e}(4,1)\times SO(3)$, $L=Spin_{e}(4,3)$.
\item[(9)] $G=SO_{e}(4,3)$, $H=SO_{e}(4,1)\times SO(2)$, $L=G_{2(2)}$.
\item[(10)] $G=SO(8,\bbC)$, $H=SO_{e}(1,7)$, $L=Spin(7,\bbC)$.
\item[(11)] $G=SO(8,\bbC)$, $H=SO(7,\bbC)$, $L=Spin_{e}(1,7)$.
\end{itemize}
\end{exes}
We now recall some features of transitive triples described in \cite{MO21}. Fix a transitive triple $(G,H,L)$. Let ${\mathfrak q}_{{\mathfrak l}}$ be the orthogonal of 
${\mathfrak l}\cap {\mathfrak h}$ in ${\mathfrak l}$ with respect to the Killing form:
\begin{equation*}
{\mathfrak l}=({\mathfrak l}\cap {\mathfrak h})\oplus {\mathfrak q}_{{\mathfrak l}}\;\text{ with }\;{\mathfrak q}_{{\mathfrak l}}=({\mathfrak l}\cap {\mathfrak h})^{\perp_{{\mathfrak l}}}\;\text{ and }\;[{\mathfrak l}\cap {\mathfrak h},{\mathfrak q}_{{\mathfrak l}}]\subset{\mathfrak q}_{{\mathfrak l}}.
\end{equation*}
Note that as ${\mathfrak l}\cap {\mathfrak h}$-modules, one has the decomposition:
\begin{equation}\label{decompoq}
{\mathfrak q}_{{\mathfrak l}}={\mathfrak q}_{{\mathfrak l}}^\prime\oplus({\mathfrak l}\cap{\mathfrak s})
\end{equation}
where ${\mathfrak q}_{{\mathfrak l}}^\prime$ is the orthogonal of ${\mathfrak l}\cap {\mathfrak h}$ in ${\mathfrak l}\cap {\mathfrak k}$ with respect to the Killing form. 
Since $L$ acts transitively on $G/H$, the $L\cap H$-equivariant map
\begin{equation*}
p^-:{\mathfrak q}_{{\mathfrak l}}\stackrel{\simeq}{\longrightarrow}{\mathfrak q},\;X\mapsto\frac{X-\sigma(X)}{2}
\end{equation*}
is an isomorphism of $L\cap H$-modules. However, the map 
\begin{equation*}
p^+:{\mathfrak q}_{{\mathfrak l}}\longrightarrow{\mathfrak h},\;X\mapsto\frac{X+\sigma(X)}{2}
\end{equation*}
need not be injective or surjective. In order to turn the isomorphism $p^-$ into an isometry, we diagonalize the symmetric bilinear form 
\begin{equation*}
{\mathfrak l}\times{\mathfrak l}\ni(X,Y)\mapsto\lan \sigma(X),Y\ran 
\end{equation*}
with respect to the Killing form $\lan ,\ran $. More precisely, for a real number $\nu$, define the vector subspace ${\mathfrak l}(\nu)$ of ${\mathfrak l}$:
\begin{equation*}
{\mathfrak l}(\nu):=\big\{X\in{\mathfrak l}\mid \lan \sigma(X),Y\ran =\nu\lan X,Y\ran \;\forall Y\in{\mathfrak l}\big\}.
\end{equation*}
Then $\fl(\nu)$ is invariant under the adjoint action of $L\cap H$ and $\fl(\nu)\perp\fl(\nu^\prime)$ whenever $\nu\neq\nu^\prime$. One can check that $\nu\in[-1,1]$, with $\fl(1)=\fl\cap\fh$ and $\fl(-1)=\fl\cap\fq$. In fact, each $\fl(\nu)$ is $\theta$-stable so that 
\begin{equation*}
\displaystyle{{\mathfrak q}_{\mathfrak l}^\prime=\bigoplus_{\nu\neq 1}{\mathfrak l}(\nu)}\cap{\mathfrak k}\;\text{ and }\;
\displaystyle{{\mathfrak l}\cap{\mathfrak s}=\bigoplus_{\nu}{\mathfrak l}(\nu)}\cap{\mathfrak s}.
\end{equation*}
We will denote by $X(\nu)$ an element in ${\mathfrak l}(\nu)$. For $\nu\neq 1$, set
\begin{equation*}
d_\nu=\Big(\frac{1-\nu}{2}\Big)^{-\frac{1}{2}}
\end{equation*}
and consider the $L\cap H$-equivariant maps
\begin{eqnarray*}
&&\rho^-:{\mathfrak q}_{\mathfrak l}\rightarrow {\mathfrak q},\; X(\nu)\mapsto d_\nu\;p^-(X(\nu))\\
&&\rho^+:{\mathfrak q}_{\mathfrak l}\rightarrow {\mathfrak h},\; X(\nu)\mapsto d_\nu\;p^+(X(\nu)).
\end{eqnarray*}
{Then $\rho^-$ is an isometry, and $\rho^+$ is a partial isometry but it need not be injective or surjective. One easily checks that for all $X(\nu)$, whenever $\nu\neq 1$ \cite{MO21}:
\begin{eqnarray*}
[X(\nu),\sigma(X(\nu))]&=& 0\\
\rho^+(X(\nu))+\rho^-(X(\nu))&=& d_\nu X(\nu)\\
\lbrack\rho^+(X(\nu)),\rho^-(X(\nu))\rbrack&=& 0.
\end{eqnarray*}
We will write $\lambda$ for the compact $\nu$, i.e the $\nu$ involved in the decomposition of ${\mathfrak q}_{\mathfrak l}^\prime$, and $\mu$ for the non-compact 
$\nu$, i.e the $\nu$ involved in the decomposition of ${\mathfrak l}\cap{\mathfrak s}$.

Transitive triples split into three families, depending on the map $\rho^+$ being injective but not surjective (type I), surjective but not injective (type S) or neither (type N).
\begin{exes}\text{}\\
\begin{itemize}
\item[(1a)] $G=G^{\prime}\times G^{\prime}$, $H=\Delta (G^{\prime}\times  G^{\prime})$, $L=G^{\prime}\times\{e\}$ is of type I.
\item[(1b)] $G=G^{\prime}\times G^{\prime}$, $H=\Delta (G^{\prime}\times G^{\prime})$, $L=G^{\prime}\times K^{\prime}$ is of type S.
\item[(2)] $G=SO_{e}(2,2n)$, $H=SO_{e}(1,2n)$, $L=U(1,n)$ is of type S.
\item[(3)] $G=SO_{e}(2,2n)$, $H=U(1,n)$, $L=SO_{e}(1,2n)$ is of type S.
\item[(4)] $G=SO_{e}(4,4n)$, $H=SO_{e}(3,4n)$, $L=Sp(1,n)$ is of type I.
\item[(5)] $G=SU(2,2n)$, $H=SU(1,2n)$, $L=Sp(1,n)$ is of type S.
\item[(6)] $G=SU(2,2n)$, $H=Sp(1,n)$, $L=SU(1,2n)$ is of type S.
\item[(7)] $G=SO_{e}(8,8)$, $H=SO_{e}(7,8)$, $L=Spin_{e}(1,8)$ is of type I.
\item[(8)] $G=SO_{e}(4,4)$, $H=SO_{e}(4,1)\times SO(3)$, $L=Spin_{e}(4,3)$ is of type N.
\item[(9)] $G=SO_{e}(4,3)$, $H=SO_{e}(4,1)\times SO(2)$, $L=G_{2(2)}$ is of type N.
\item[(10)] $G=SO(8,\bbC)$, $H=SO_{e}(1,7)$, $L=Spin(7,\bbC)$ is of type N.
\item[(11)] $G=SO(8,\bbC)$, $H=SO(7,\bbC)$, $L=Spin_{e}(1,7)$ is of type I.
\end{itemize}
\end{exes}
We fix an orthonormal basis $\{Z_j(\lambda)\}$ (resp. $\{T_k(\mu)\}$) of  $\fl(\lambda)$ (resp. $\fl(\mu)$) so that $\{Z_j\}=\cup_\lambda\{Z_j(\lambda)\}$ (resp. $\{T_k\}=\cup_\mu\{T_k(\mu)\}$) is an orthonormal basis of ${\mathfrak q}_{\mathfrak l}^\prime$ (resp. ${\mathfrak l}\cap{\mathfrak s}$) and
\begin{eqnarray}\label{thebasis}
&&[Z_j(\lambda),\sigma(Z_j(\lambda))]=0,\;\;\forall j,\lambda.\nonumber\\
&&[T_k(\mu),\sigma(T_k(\mu))]=0,\;\;\forall k,\mu.
\end{eqnarray}
It is worth to mention that transitive triples of type S have the following additional features \cite{MO21}:
\begin{itemize}
\item[(a)] $L\cap K$ is $\sigma$-stable, i.e., $L\cap K/L\cap H$ is a symmetric space with respect to the involution $\sigma$.
\item[(b)] ${\mathfrak q}_{\mathfrak l}^\prime$ is irreducible as a $L\cap H$-module and $\lambda=-1$ (only one value for $\lambda$).
\item[(c)] ${\mathfrak l}\cap{\mathfrak s}$ is irreducible as a $L\cap H$-module and $\mu=0$ (only one value for $\mu$).
\end{itemize}
%
\section{Clifford algebras and spin modules}\label{Spin}

Let $(G,H,L)$ be a transitive triple, and let as before $\frq$ denote the orthogonal complement of $\frh$ is $\frg$ with respect to the form $\lara$. Then $\lara$ is nondegenerate on both $\frh$ and $\frq$.

Let $C(\frq)$ be the Clifford algebra of $\frq$ with respect to $\lara$. In other words, it is the associative algebra with unit, generated by $\frq$, with relations
\[
XY+YX=\lan X,Y\ran,\qquad X,Y\in\frq.
\]
It is well known that the category of complex $C(\frq)$-modules is semisimple, with only one irreducible module if $\dim\frq$ is even, and two irreducible modules if $\dim\frq$ is odd. These modules are called spin modules and they can be constructed as follows. Let $\frq^+_\bbC$ and $\frq^-_\bbC$ be maximal isotropic subspaces of $\frq_\bbC$, nondegenerately paired by $\lara$. Then 
\begin{gather*}
\frq_\bbC=\frq^+_\bbC\oplus\frq^-_\bbC\qquad\qquad\ \ \text{if $\dim\frq$ is even};\\
\frq_\bbC=\frq^+_\bbC\oplus\frq^-_\bbC\oplus\bbC Z\qquad\text{if $\dim\frq$ is odd},
\end{gather*}
where in the odd case $Z\in\frq_\bbC$ is a vector orthogonal to $\frq^+_\bbC\oplus\frq^-_\bbC$ such that $\lan Z,Z\ran=1$.
We define
\[
S_\frq=\twedge\frq^+_\bbC,
\]
with elements of $\frq^+_\bbC$ acting by wedging and elements of $\frq^-_\bbC$ acing by contracting. If $\dim\frq$ is even, this determines a $C(\frq)$-module structure on the spin module 
$S_\frq$, and if $\dim\frq$ is odd we still need to define the action of $Z$. There are two choices: $Z$ can act by $1/\sqrt{2}$ on $\twedge^{\even}\frq^+_\bbC$ and by $-1/\sqrt{2}$ on $\twedge^{\odd}\frq^+_\bbC$, or by $-1/\sqrt{2}$ on $\twedge^{\even}\frq^+_\bbC$ and by $1/\sqrt{2}$ on $\twedge^{\odd}\frq^+_\bbC$. We fix one of these choices.
We denote by 
\begin{equation}\label{defcliffmult}
\gamma_\fq:C(\fq)\to\End(S_\fq)
\end{equation}
the action map for the $C(\frq)$-module $S_\frq$.

Besides the obvious embedding of $\frq$ into $C(\frq)$, there is also an embedding of $\frso(\frq)$ into $C(\frq)$ induced by the skew symmetrization (Chevalley isomorphism):
\begin{equation*}
\fs\fo(\fq)\simeq\Lambda^2\fq\overset{j}{\hookrightarrow} C(\fq),\qquad j(X\wedge Y)=\frac{1}{2}(XY-YX).
\end{equation*}
The image of this map is equal to the Lie subalgebra $C^2(\fq)$ of $C(\frq)$, consisting of ``pure degree $2$" elements. Its main property is that the natural action of $\frso(\frq)$ on $\frq$ corresponds to the action of $C^2(\frq)$ on $\frq\subset C(\frq)$ by Clifford algebra commutators.

Composing the above map with the adjoint action map $\frh\to\frso(\frq)$, we get a map
\begin{equation}\label{alphah}
\alpha_{\frg,\fh}=\alpha_\fh:\fh\to C(\fq)
\end{equation}
such that
\eq\label{cliffbracket}
[\alpha_\fh(X),Y]_{C(\fq)}=[X,Y]_\frg,\qquad X\in\frh,Y\in\frq.
\eeq
Moreover, $\alpha_\frh$ is a Lie algebra morphism, i.e., 
\[
[\alpha_\fh(X),\alpha_\fh(Z)]_{C(\fq)}=[X,Z]_\frg,\qquad X,Z\in\frh.
\]
In the last two formulas, the subscript of a bracket denotes the Lie algebra in which the bracket is taken. 

The following lemma describes $\alpha_\frh$ more explicitly. The proof can be found for example in \cite{hpbook}, Section 2.3.3 (but note that the conventions there are different, so a factor shows up in the formula).

\begin{lem}
\label{alphah2}
Let $\{e_i\}$ be an orthonormal basis of $\fq$, and let
$\eps_i=\lan e_i,e_i\ran\in\{\pm 1\}$. Then for any $X\in\fh$,
\[
\alpha_\frh(X)=-\sum_{i<j}\eps_i\eps_j \lan [X,e_i],e_j\ran e_ie_j=-\frac{1}{2}\sum_{i,j}\eps_i\eps_j \lan [X,e_i],e_j\ran e_ie_j.
\]
\end{lem}
Using the map $\alpha_\frh$, we can view the spin module $S_\frq$ as an $\frh$-module, with $X\in\frh$ acting as $\alpha_\frh(X)\in C(\frq)$.

Recall that $\rho^-:\frq_\frl\to\frq$ is an isometry, and therefore $\rho^-$ extends to an isomorphism of 
the Clifford algebras $C({\mathfrak q}_{\mathfrak l})$ and $C({\mathfrak q})$, where $\rho^-:C({\mathfrak q}_{\mathfrak l})\to C({\mathfrak q})$ is defined by
\[
\rho^-(Y_1\dots Y_k)=\rho^-(Y_1)\dots \rho^-(Y_k),\qquad Y_1,\dots,Y_k\in \frq_\frl.
\]
It follows that also the spin modules $S_{\frq_\frl}$ and $S_\frq$ are isomorphic as vector spaces. On the other hand, $S_{\frq_\frl}$ has an action of $\frl\cap\frh$, while $S_\frq$ has an action of $\frh$, and hence also of $\frl\cap\frh$ by restriction.

\begin{lem}
\label{l cap h actions}
With notation as above,
\[
\rho^-(\alpha_{\frl\cap\frh}(X))=\alpha_\frh(X),\qquad X\in\frl\cap\frh.
\]
In particular, the isomorphism of $S_{\frq_\frl}$ and $S_\frq$  induced by $\rho^-$ is $\frl\cap\frh$-equivariant.
\end{lem}
\pf Recall that $\rho^-:\frq_\frl\to\frq$ is an $\frl\cap\frh$-equivariant isometry. This implies that  for any $X\in\frl\cap\frh$ and $Y,Z\in\frq_\frl$,
\eq
\label{l cap h rho-}
\lan[X,\rho^-(Y)],\rho^-(Z)\ran=\lan\rho^-([X,Y]),\rho^-(Z)\ran=\lan[X,Y],Z\ran.
\eeq
Let now $\{e_i\}$ be an orthonormal basis of $\frq_l$ with $\lan e_i,e_i\ran=\eps_i=\pm 1$. By Lemma \ref{alphah2}, if $X\in\frl\cap\frh$, then
\begin{gather*}
\alpha_{\frl\cap\frh}(X)=-\sum_{i<j}\eps_i\eps_j \lan [X,e_i],e_j\ran e_ie_j;\\
\alpha_\frh(X)=-\sum_{i<j}\eps_i\eps_j \lan [X,\rho^-(e_i)],\rho^-(e_j)\ran \rho^-(e_i)\rho^-(e_j).
\end{gather*}
So the lemma follows from \eqref{l cap h rho-}.
\epf
We also consider the spin modules $S_{{\mathfrak q}_{\mathfrak l}^\prime}$ and $S_{{\mathfrak l}\cap{\mathfrak s}}$ for the pairs $({\mathfrak l}\cap{\mathfrak k},{\mathfrak l}\cap{\mathfrak h})$ and $({\mathfrak l},{\mathfrak l}\cap{\mathfrak k})$ respectively. Recall that by (\ref{decompoq})
\[
\frq_\frl=\frl\cap\frs\oplus\frq_\frl'.
\]
This implies that 
\[
C(\frq_\frl)= C(\frl\cap\frs)\,\bar\otimes\, C(\frq_\frl'),
\]
where $\bar\otimes$ denotes the graded tensor product of superalgebras, i.e.,
\[
(a\bar\otimes b)(c\bar\otimes d)=(-1)^{{\scriptsize\deg}\, b\,{\scriptsize\deg}\, c}ac\,\bar\otimes\, bd,
\]
with $\deg$ being 0 for even elements and 1 for odd elements.

To get an analogous decomposition of the spin module, we first note that if $(\frl\cap\frs)^+_\bbC$ and $(\frq_\frl')^+_\bbC$ are maximal isotropic subspaces of $(\frl\cap\frs)_\bbC$ respectively $(\frq_\frl')_\bbC$, then $(\frl\cap\frs)^+_\bbC\oplus(\frq_\frl')^+_\bbC$ is a maximal isotropic subspace of $(\frq_\frl)_\bbC$, unless $\dim(\frl\cap\frs)$ and $\dim(\frq_\frl')$ are both odd. (In case they are both odd, there is an isotropic vector outside of $(\frl\cap\frs)^+_\bbC\oplus(\frq_\frl')^+_\bbC$.) To simplify matters, we make the following assumption. 
\begin{assump}\label{spinassum}
The dimension of $\frl\cap\frs$ is even.
\end{assump}
Since the triples of type $S$ other than 
$(G,H,L)=(G^\prime\times G^\prime,\Delta(G^\prime\times G^\prime),G^\prime\times K^\prime)$ all have $L$ of equal rank,
Assumption \ref{spinassum} is satisfied for all of them. For the triples $(G^\prime\times G^\prime,\Delta(G^\prime\times G^\prime),G^\prime\times K^\prime)$, we have $\frl\cap\frs=\frs'\times 0$, so Assumption \ref{spinassum} is equivalent to $\dim\frs'$ being even. This is certainly true if $G'$ is of equal rank, and if $G'$ is of unequal rank it is still often true. For example, if $G'=SL(n,\bbR)$, then $\dim\frs'=\frac{n(n+1)}{2}-1$, which is even if $n$ is congruent to 1 or 2 modulo 4 (and odd otherwise).

The above discussion about isotropic subspaces implies that under Assumption \ref{spinassum} we have
$S_{\mathfrak q_\frl}= S_{{\mathfrak l}\cap{\mathfrak s}} \otimes S_{{\mathfrak q}_{\mathfrak l}^\prime}$
as vector spaces. The isomorphism 
$S_{{\mathfrak l}\cap{\mathfrak s}} \otimes S_{{\mathfrak q}_{\mathfrak l}^\prime}\to S_{{\mathfrak q}_{\mathfrak l}}$ can be realized as the exterior algebra multiplication map.

This vector space isomorphism is also an isomorphism of $C(\frq_\frl)=C(\frl\cap\frs)\,\bar\otimes\, C(\frq_\frl')$-modules. To see this, we first note that since $\dim\frl\cap\frs$ is even, the $C(\frl\cap\frs)$-module $S_{\frl\cap\frs}$ is $\bbZ_2$-graded, with
\[
S_{\frl\cap\frs}^0=\twedge^{\text{even}}(\frl\cap\frs)^+_\bbC;\qquad S_{\frl\cap\frs}^1=\twedge^{\text{odd}}(\frl\cap\frs)^+_\bbC.
\] 
We can now define an action of $C(\frq_\frl)=C(\frl\cap\frs)\,\bar\otimes\, C(\frq_\frl')$ on $S_{{\mathfrak l}\cap{\mathfrak s}} \otimes S_{{\mathfrak q}_{\mathfrak l}^\prime}$ by
\eq
\label{spindec action}
(c\otimes d)\cdot (s\otimes t)=(-1)^{{\scriptsize\deg}\, d\,{\scriptsize\deg}\, s}cs\otimes dt.
\eeq
To see that this is a well defined algebra action, we check
\eq
\label{spin dec wd}
(a\otimes b)\cdot[(c\otimes d)\cdot(s\otimes t)]=[(a\otimes b)(c\otimes d)]\cdot(s\otimes t),
\eeq
for any $a,c\in C(\frl\cap\frs)$, $b,d\in C(\frq_\frl')$, $s\in S_{\frl\cap\frs}$ and $t\in S_{\frq_\frl'}$ such that $b,c,d$ and $s$ are homogeneous.

Up to sign, both sides of \eqref{spin dec wd} are equal to 
$acs\otimes bdt$. The sign on the left side is $-1$ to the power of $\deg d\deg s + \deg b\deg(cs)$, while the sign on the right side is $-1$ to the power of $\deg b\deg c + \deg (bd)\deg s$. Both exponents are equal to 
\[
\deg b\deg c+\deg b\deg s+\deg d\deg s,
\]
so \eqref{spin dec wd} is true and the action \eqref{spindec action} is well defined. Moreover, we have
\begin{lem}
\label{spindec cliff}
The exterior algebra multiplication map
\[
m:S_{{\mathfrak l}\cap{\mathfrak s}} \otimes S_{{\mathfrak q}_{\mathfrak l}^\prime}\to S_{{\mathfrak q}_{\mathfrak l}}
\]
is an isomorphism of $C(\frq_\frl)=C(\frl\cap\frs)\,\bar\otimes\, C(\frq_\frl')$-modules, where the action of $C(\frq_\frl)$ on $S_{{\mathfrak q}_{\mathfrak l}}$ is the usual one, and the action on $S_{{\mathfrak l}\cap{\mathfrak s}} \otimes S_{{\mathfrak q}_{\mathfrak l}^\prime}$ is the one defined by \eqref{spindec action}.
\end{lem}
\pf We already know that $m$ is a vector space isomorphism, so we only need to check it respects the actions. This follows immediately from the definitions.
\epf

It is now not difficult to see that we also have 
\begin{equation}\label{isomspin}
S_{\frq_\frl}= S_{{\frl}\cap{\frs}} \otimes 
S_{\frq_\frl'}
\end{equation}
on the level of $\frl\cap\frh$-modules, where 
$S_{\frq_\frl}$ is viewed as an $\frl\cap\frh$-module through the map $\alpha_{\frl,\frl\cap\frh}:\frl\cap\frh\to C(\frq_\frl)$, $S_{\frl\cap\frs}$ is viewed as an $\frl\cap\frh$-module through the restriction of the map $\alpha_{\frl,\frl\cap\frk}:\frl\cap\frk\to C(\frl\cap\frs)$, and $S_{\frq_\frl'}$ is viewed as an $\frl\cap\frh$-module through the map $\alpha_{\frl\cap\frk,\frl\cap\frh}:\frl\cap\frh\to C(\frq_\frl')$. Indeed, we have
\begin{lem}
\label{spindec}
Let $(G,H,L)$ be a triple of type $S$ satisfying Assumption \ref{spinassum}. Then for any $X\in\frl\cap\frh$, $\alpha_{\frl,\frl\cap\frh}(X)$ decomposes 
under $C(\frq_\frl)=C(\frl\cap\frs)\,\bar\otimes\, C(\frq_\frl')$ as
\[
\alpha_{\frl,\frl\cap\frh}(X)=\alpha_{\frl,\frl\cap\frk}(X)\otimes 1+1\otimes\alpha_{\frl\cap\frk,\frl\cap\frh}(X).
\]
In particular, the isomorphism of Lemma \ref{spindec cliff} is an isomorphism of $\frl\cap\frh$-modules.
\end{lem}
\pf
Let $\{Z_i\}$ respectively $T_r$ be orthonormal bases of $\frq_\frl'$ respectively $\frl\cap\frs$. Since $\lan Z_i,Z_i\ran=-1$ and $\lan T_k,T_k\ran=1$ for any $i$ and $k$, Lemma \ref{alphah2} implies that for any $X\in\frl\cap\frh$, 
\eq
\label{spin dec 2}
\alpha_{\frl,\frl\cap\frh}(X)=-\sum_{i<j}\lan [X,Z_i],Z_j\ran Z_iZ_j +\sum_{i,r}\lan [X,Z_i],T_r\ran Z_iT_r - \sum_{r<s}\lan [X,T_r],T_s\ran T_rT_s.
\eeq 
Since $[X,Z_i]\in\frl\cap\frk\perp\frl\cap\frs$, the second sum in \eqref{spin dec 2} is 0, i.e., 
\eq
\label{spin dec 3}
\alpha_{\frl,\frl\cap\frh}(X)=-\sum_{i<j}\lan [X,Z_i],Z_j\ran Z_iZ_j - \sum_{r<s}\lan [X,T_r],T_s\ran T_rT_s.
\eeq 
Since the first sum in \eqref{spin dec 3} is equal to 
$\alpha_{\frl\cap\frk,\frl\cap\frh}(X)=1\otimes\alpha_{\frl\cap\frk,\frl\cap\frh}(X)$ by Lemma \ref{alphah2}, 
and the second sum is equal to $\alpha_{\frl,\frl\cap\frk}(X)=\alpha_{\frl,\frl\cap\frk}(X)\otimes 1$, the lemma follows.
\epf


\section{Twisted Dirac operators}\label{Dirac}

As before, let $G$ be a non-compact connected semisimple real Lie group with finite center and with Lie algebra ${\mathfrak g}$, $K$ a maximal compact subgroup of $G$ with respect to a Cartan involution $\theta$ and $H$ a non-compact connected semisimple closed subgroup of $G$ with Lie algebra ${\mathfrak h}$ such that $G/H$ is a symmetric space with respect to an involution $\sigma$ commuting with $\theta$.

Let $U(\fg_{\mathbb C})$ (resp. $U(\fh_{\mathbb C})$) be the enveloping algebra of the complexification of $\fg$ (resp. $\fh$). For $i=1,2$, let 
\begin{equation*}
\tau_{i}:H\rightarrow GL(F_{i})
\end{equation*}
be finite dimensional smooth complex representations of $H$. Write $\Hom(F_{1},F_{2})$ for the vector space of complex homomorphisms from $F_{1}$ to $F_{2}$. The enveloping algebra $U(\fg_{\mathbb C})$ is both a right $U(\fh_{\mathbb C})$-module and an $H$-module: 
\begin{eqnarray*}
X\cdot A&=&AX\\
h\cdot A&=&\Ad(h)A
\end{eqnarray*} 
for all $h\in H$, $X\in U(\fh_{\mathbb C})$ and $A\in U(\fg_{\mathbb C})$. Write $X\mapsto X^{0}$ for the anti-automorphism of $U(\fh_{\mathbb C})$ defined by 
\begin{equation*}
(X_{1}\cdots X_{n})^{0}=(-1)^{n}X_{n}\cdots X_{1}\quad\text{ for }X_{j}\in\fh.
\end{equation*}
In particular, the vector space $\Hom(F_{1},F_{2})$ is equipped with a structure of $U(\fh_{\mathbb C})$-module: 
\begin{equation}\label{action0}
X\cdot T=T\circ d\tau_{1}(X^{0}),\qquad X\in U(\fh_{\mathbb C}), \;T\in\Hom(F_{1},F_{2})
\end{equation}
where $d\tau_{1}$ denotes the differential of $\tau_{1}$ extended naturally to $U(\fh_{\mathbb C})$. Let ${\Cal F}_{i}\rightarrow G/H$ be the homogeneous vector bundle over $G/H$ induced by the $H$-module $F_{i}$ and write $C^{\infty}(G/H,{\Cal F}_{i})$ for the space of smooth sections on which the group $G$ acts by left translations. As a $G$-module, $C^{\infty}(G/H,{\Cal F}_{i})$ is isomorphic to the space of $H$-invariant vectors $\big(C^{\infty}(G)\otimes F_{i}\big)^{H}$. Here the $G$-action on $C^{\infty}(G)\otimes F_{i}$ is given by left translations on $C^{\infty}(G)$ and trivial on $F_i$, while the $H$-action is given by the right translations on $C^{\infty}(G)$ and by $\tau_{i}$ on $F_{i}$. Let ${\mathcal D}_{G/H}({\Cal F}_{1},{\Cal F}_{2})$ be the vector space of left-invariant differential operators $C^{\infty}(G/H,{\Cal F}_{1})\rightarrow C^{\infty}(G/H,{\Cal F}_{2})$. One has the following isomorphism:
\begin{equation*}
{\mathcal D}_{G/H}({\Cal F}_{1},{\Cal F}_{2})\simeq\Big\{U(\fg_{\mathbb C})\otimes_{U(\fh_{\mathbb C})}\Hom(F_{1},F_{2})\Big\}^{H}
\end{equation*}
where the $H$-action on $U(\fg_{\mathbb C})\otimes_{U(\fh_{\mathbb C})}\Hom(F_{1},F_{2})$ is given by 
\begin{equation}\label{action}
h\cdot (A\otimes T)=(\Ad(h)A)\otimes\tau_{2}(h)\circ T\circ\tau_{1}(h)^{-1}
\end{equation}
When $F_1=F_2$, ${\mathcal D}_{G/H}({\Cal F}_{1},{\Cal F}_{2})$ is an algebra and in the case when $F_{1}=F_{2}=\bbC$, one has an isomorphism of algebras
\beu
{\mathcal D}_{G/H}(\bbC,\bbC)\simeq U(\fg_{\mathbb C})^{H}/U(\fg_{\mathbb C})^{H}\cap U(\fg_{\mathbb C})\fh_{\mathbb C}
\end{equation*}
where ${\mathcal D}_{G/H}(\bbC,\bbC)$ coincides with the commutative algebra ${\mathcal D}(G/H)$ of left-invariant differential operators acting on smooth functions on $G/H$. 

Let $(\beta,E)$ be a finite dimensional representation of $\fh$ such that the tensor product $S_{\fq}\otimes E$ lifts to a representation $(\tau, F)$ of the group $H$. There is an associated smooth homogeneous vector bundle over $G/H$, which we denote by $\Cal S_\fq\otimes\Cal E$, whose space of smooth sections is 
\begin{eqnarray*}
 && C^\infty(G/H,\Cal S_\fq\otimes\Cal E)\\
   & \simeq& \Big\{C^{\infty}(G)\otimes(S_{\fq}\otimes E)\Big\}^{H}  \\
    &\simeq& \{f: G\to S_\fq\otimes E\, |\, f\text{ is smooth and }f(gh)=\tau(h)^{-1}(f(g)),\; \forall h\in H\} .
  \end{eqnarray*}
For $X\in{\mathfrak q}$ and $\phi\in C^\infty(G)$, define the right differential of $\phi$ along $X$ as follows:
\begin{equation*}
(R(X))\varphi)(g)=\frac{d}{dt}\varphi(g\exp(tX))_{\mid_{t=0}}.
\end{equation*}
Pick an orthonormal basis $\{X_j\}$ of $\fq$ and consider the operator 
\begin{equation}\label{defdirac}
\widehat{D}_{G/H}(E):C^\infty(G/H,\Cal S_\fq\otimes\Cal E)\longrightarrow C^\infty(G/H,\Cal S_\fq\otimes\Cal E)
\end{equation}
defined by
\begin{equation*}
\widehat{D}_{G/H}(E)=\sum_j\lan X_j,X_j\ran R(X_j)\otimes(\gamma_\fq(X_j)\otimes 1).
\end{equation*}
Note that $\widehat{D}_{G/H}(E)$ is independent of the basis $\{X_j\}$. One checks that $\widehat{D}_{G/H}(E)$ belongs to the space $\Big\{U(\fg_{\mathbb C})\otimes_{U(\fh_{\mathbb C})}\Hom(S_\fq\otimes E,S_\fq\otimes E)\Big\}^{H}$ and therefore defines a $G$-invariant differential operator acting on $C^\infty(G/H,\Cal S_\fq\otimes \Cal E)$. 
$\widehat{D}_{G/H}(E)$ is known as the {\it twisted geometric Dirac operator} associated with $E$. 

There is an algebraic analog of the geometric Dirac operator. Namely, attached to a $\fg$-module $(\pi,V)$ there is the Dirac operator
\begin{equation*}
D_{\fg,\fh}(V):V\otimes S_\fq\rightarrow V\otimes S_\fq
\end{equation*}
with 
\begin{equation}\label{algdirac}
D_{\fg,\fh}(V):=\sum_j\lan X_j,X_j\ran \pi(X_j)\otimes\gamma_\fq(X_j).
\end{equation}
The operator $D_{\fg,\fh}(V)$ is independent of the choice of basis $\{X_j\}$ and is $\frh$-invariant.

\section{Embedding of Dirac operators}\label{embedding}

Fix a transitive triple $(G,H,L)$. Recall that the spin modules $S_{{\mathfrak q}_{\mathfrak l}}$ and $S_{\mathfrak q}$ for ${\mathfrak l}\cap{\mathfrak h}$ and ${\mathfrak h}$ respectively are isomorphic as ${\mathfrak l}\cap{\mathfrak h}$-modules. Let $(\beta,E)$ be a finite dimensional representation of $\fh$ such that the tensor product $S_{\fq}\otimes E$ lifts to a representation $(\tau, F)$ of the group $H$. There is an associated smooth homogeneous vector bundle over both $G/H$ and $L/L\cap H$, which we denote by $\Cal S_\fq\otimes\Cal E\rightarrow G/H$ and $\Cal S_{\fq_\frl}\otimes\Cal E\rightarrow L/L\cap H$ respectively. The spaces of smooth sections are, as above, denoted respectively by $C^\infty(G/H,\Cal S_\fq\otimes\Cal E)$ and 
$C^\infty(L/L\cap H,\Cal S_{\fq_\fl}\otimes\Cal E)$. The transitive action of $L$ on $G/H$ implies that restriction to $L$ is an isomorphism of $L$-modules:
\begin{equation}\label{mapj}
\jmath:C^\infty(G/H,\Cal S_\fq\otimes\Cal E)\stackrel{\simeq}{\longrightarrow} C^\infty(L/L\cap H,\Cal S_{\fq_\fl}\otimes\Cal E),\quad f\mapsto f_{\mid_L}.
\end{equation}

On the other hand, if $U(\fl_{\mathbb C})$ denotes the enveloping algebra of the complexification $\fl$, the transitive action of $L$ on $G/H$ and the Poincar\'e-Birkhoff-Witt theorem induce an isomorphism of algebras 
\begin{equation*}
U({\mathfrak g}_{\mathbb C})\simeq U({\mathfrak h}_{\mathbb C})\otimes_{U({\mathfrak l}_{\mathbb C}\cap{\mathfrak h}_{\mathbb C})}U({\mathfrak l}_{\mathbb C}).
\end{equation*}
One deduces the following $L\cap H$-equivariant embedding of invariant differential operators:
\begin{gather}\label{embedop}
\imath:{\mathcal D}_{G/H}({\mathcal S}_{\mathfrak q}\otimes {\mathcal E},{\mathcal S}_{\mathfrak q}\otimes {\mathcal E})\hookrightarrow {\mathcal D}_{L/L\cap H}({\mathcal S}_{{\mathfrak q}_{\mathfrak l}}\otimes {\mathcal E},{\mathcal S}_{{\mathfrak q}_{\mathfrak l}}\otimes {\mathcal E})\nonumber\\
\imath(D)(\jmath(f))=\jmath(D(f)).
\end{gather}
Moreover, recall that under Assumption \ref{spinassum}, the spin modules $S_{\fq_\fl^\prime}$ for the pair $(\fl\cap\fk,\fl\cap\fh)$ and $S_{\fl\cap\fs}$ for the pair 
$(\fl,\fl\cap\fk)$ satisfy the isomorphism (\ref{isomspin}) of $\fl\cap\fh$-modules: $S_{\fq_\fl}\simeq S_{\fl\cap\fs}\otimes S_{\fq_\fl^\prime}$. We get the following isomorphism of $L$-modules:
\begin{gather}
\Phi:C^\infty(L/L\cap H,\Cal S_{\fq_\fl}\otimes\Cal E)\stackrel{\simeq}{\longrightarrow} C^\infty(L/L\cap K,\Cal S_{\fl\cap\fs}\otimes C^\infty(L\cap K/L\cap H,\Cal S_{\fq_\fl^\prime}\otimes\Cal E))\nonumber\\
\label{isomPhi}
\Phi(f)(l)(k)=(k\otimes 1)(f(lk))\;\;\forall l\in L,\; k\in L\cap K
\end{gather}
where $k$ acts on the first factor of $S_{\fl\cap\fs}\otimes (S_{\fq_\fl^\prime}\otimes E)$. Define the $L\cap K$-module $\widetilde{E}$ and $\widetilde{\mathcal E}$ the corresponding homogeneous vector bundle over $L\cap K/L\cap H$:
\begin{equation}\label{tildeE}
\widetilde{E}:=C^\infty(L\cap K/L\cap H,\Cal S_{\fq_\fl^\prime}\otimes\Cal E)\;\text{ and }\;\widetilde{\mathcal E}\rightarrow L\cap K/L\cap H.
\end{equation}

Next, as in (\ref{defdirac}), one can define the operator 
\begin{equation}\label{noncubicDiracLH}
\widehat{D}_{L/L\cap H}(E):C^\infty(L/L\cap H,{\mathcal S}_{\fq_\fl}\otimes{\mathcal E})\rightarrow C^\infty(L/L\cap H,{\mathcal S}_{\fq_\fl}\otimes{\mathcal E}).
\end{equation}
Let $c_{\fl,\fl\cap\fh}$ be the degree three element in $C(\fq_\fl)$ defined as the image under the Chevalley isomorphism of the $3$-form on $\fq_\fl$ given by 
\begin{equation*}
{\mathfrak q}_\fl\times\fq_\fl\times\fq_\fl\ni(X,Y,Z)\mapsto\IP{X}{\lbrack Y,Z\rbrack}.
\end{equation*}
If $\{e_j\}$ is an orthonormal basis of ${\mathfrak q}_\fl$ with $\lan e_j,e_j\ran=\eps_j=\pm 1$, then 
\begin{equation}\label{cubicterm}
c_{\fl,\fl\cap\fh}=\sum_{i<j<k}a_i a_j a_k\lan [e_i,e_j],e_k\ran e_i e_j e_k=\frac{1}{6}\sum_{i,j,k}a_i a_j a_k\lan [e_i,e_j],e_k\ran e_i e_j e_k.
\end{equation}
Note that both $\widehat{D}_{L/L\cap H}(E)$ and $c_{\fl,\fl\cap\fh}$ belong to the space
\[
\Big\{U(\fl_{\mathbb C})\otimes_{U(\fl_{\mathbb C}\cap\fh_{\mathbb C})}
\Hom(S_{\fq_\fl}\otimes E,S_{\fq_\fl}\otimes E)\Big\}^{L\cap H}.
\] 
The $L$-invariant differential operator $\widehat{D}_{L/L\cap H}(E)$ is the non-cubic geometric Dirac operator and 
\begin{equation}\label{cubicDiracLH}
D_{L/L\cap H}(E):=\widehat{D}_{L/L\cap H}(E)-1\otimes\gamma_{\fq_\fl}(c_{\fl,\fl\cap\fh})\otimes 1
\end{equation}
is the cubic geometric Dirac operator. 

Similarly, following (\ref{algdirac}), one can define respectively the cubic and non-cubic algebraic Dirac operator $D_{\fl,\fl\cap\fh}(E)$ and $\widehat{D}_{\fl,\fl\cap\fh}(E)$ on $S_{\fq_\fl}\otimes E$:
\begin{equation*}
D_{\fl,\fl\cap\fh}(E)=\widehat{D}_{\fl,\fl\cap\fh}(E)-1\otimes\gamma_{\fq_\fl}({c}_{\fl,\fl\cap\fh}).
\end{equation*}
The algebraic Dirac operator for $(\fh,\fh\cap\fk)$ is defined analogously. Note that since $G/H$ (resp. $K/K\cap H$) is a symmetric space, the cubic term ${c}_{\fg,\fh}$ (resp. ${c}_{\fh,\fh\cap\fk}$) for the pair $(\fg,\fh)$ (resp. $(\fh,\fh\cap\fk)$) vanishes, so that cubic and non-cubic Dirac operators coincide. In particular, one may write equivalently $D_{G/H}(E)$ or $\widehat{D}_{G/H}(E)$. Using the isomorphism $\Phi$, the operator $D_{L/L\cap H}(E)$ may be pushed over to the right side of \eqref{isomPhi} to an $L$-invariant differential operator:
\begin{equation}\label{pushD}
\widetilde{D}_{L/L\cap H}(E)(\Phi(f))=\Phi(D_{L/L\cap H}(E)(f)).
\end{equation}
We will omit the tilde and still write $D_{L/L\cap H}(E)$ for the pushed operator $\widetilde{D}_{L/L\cap H}(E)$ when there is no confusion.
\begin{thm}\label{embeddingthm} For transitive triples of type $S$, one has:
\begin{multline*}
\imath(D_{G/H}(E))=\sqrt{2}\;D_{L/L\cap H}(E)+(1-\sqrt{2})D_{L\cap K/L\cap H}(E)+\\
(\sqrt{2}-2)\otimes \gamma_{\frq_\frl}(c_{\fl,\fl\cap\fh})\otimes 1+1\otimes D_{\fh,\fh\cap\fk}(E).
\end{multline*}
Another way to express $\imath(D_{G/H}(E))$ is
\[
\imath(D_{G/H}(E))=\sqrt{2}\;D_{L/L\cap K}(\widetilde{E})+D_{L\cap K/L\cap H}(E)-2\otimes \gamma_{\frq_\frl}(c_{\fl,\fl\cap\fh})\otimes 1+1\otimes D_{\fh,\fh\cap\fk}(E).
\]
The actions of $D_{L/L\cap H}(E)$ and $c_{\fl,\fl\cap\fh}$ were defined above, while $D_{L\cap K/L\cap H}(E)$ acts on $\widetilde{E}=C^\infty(L\cap K/L\cap H,{\mathcal S}_{\fq_\fl}\otimes{\mathcal E})$, $D_{L/L\cap K}(\widetilde{E})$ acts on $C^\infty(L/L\cap K,{\mathcal S}_{\fl\cap\fs}\otimes\widetilde{\mathcal E})$ and $D_{\fh,\fh\cap\fk}(E)$ acts on $S_{\fh\cap\fs}\otimes E\simeq S_{\fl\cap\fs}\otimes E$.
\end{thm}

\pf
Recall that for triples of type $S$ the only $\lambda$ is $-1$ with $d_\lambda=1$ and the only $\mu$ is $0$ with $d_\mu=\sqrt{2}$. We denote by $\{Z_j\}$ and $\{T_k\}$ 
the orthonormal bases  of ${\mathfrak q}_\fl'$ and $\fl\cap\fs$ defined in (\ref{thebasis}). The Dirac operator $D_{G/H}(E)$ is then given by
\begin{equation*}
D_{G/H}(E)=-\sum_j\rho^-(Z_j)\otimes\gamma_{\fq}(\rho^-(Z_j))\otimes 1+\sum_k\rho^-(T_k)\otimes\gamma_{\fq}(\rho^-(T_k))\otimes 1.
\end{equation*}
It is an $H$-invariant element of 
\[
U(\fg_{\mathbb C})\otimes_{U(\fh_{\mathbb C})}\text{End}(S_\fq\otimes E)=U(\fg_{\mathbb C})\otimes_{U(\fh_{\mathbb C})}\text{End}(S_\fq)\otimes \text{End}(E).
\]

We know from Section \ref{triples} that 
\begin{gather*}
\rho^-(Z_j)=Z_j,\qquad \rho^+(Z_j)=0;\\
\rho^-(T_k)=\sqrt{2}\,T_k -\rho^+(T_k).
\end{gather*} 
Using (\ref{action0}) and remembering that the $\frh$-action on $S_\frq$ is through the map $\alpha_\frh$ of (\ref{alphah}), we get

\begin{multline*}
D_{G/H}(E)=\underbrace{-\sum_j  Z_j\otimes\gamma_{\fq}(\rho^-(Z_j))\otimes 1+\sqrt{2}\sum_k T_k\otimes\gamma_{\fq}(\rho^-(T_k))\otimes 1}_{D_1}\\
+1\otimes\gamma_{\fq}\Big(\underbrace{\sum_k \rho^-(T_k)\alpha_\fh(\rho^+(T_k))}_{D_2}\Big)\otimes 1
+1\otimes\underbrace{\Big(
\sum_k \gamma_{\fq}(\rho^-(T_k))\otimes\beta(\rho^+(T_k))\Big)}_{D_3},
\end{multline*}
where $\beta(\rho^+(T_k))\in\End(E)$ is the action of $\rho^+(T_k)\in\frh$ on the $\frh$-module $E$.

It follows that
\[
\imath(D_{G/H}(E))=\imath(D_1)+\imath(1\otimes \gamma_\frq(D_2)\otimes 1)+\imath(1\otimes D_3).
\]
Recall that the embedding $\imath$ includes using
the isometry $\rho^-:\frq_\frl\to\frq$ to identify the Clifford algebra and spin module for 
$\frq_{\frl}$ with the Clifford algebra and spin module for $\frq$. Effectively this means that in the Clifford algebra factor of the above expressions we replace $\gamma_{\fq}(\rho^-(Z_j))$ respectively $\gamma_{\fq}(\rho^-(T_k))$ with $\gamma_{\fq_\frl}(Z_j)$ respectively $\gamma_{\fq_\frl}(T_k)$.

In particular, 
\[
\imath(D_1)=-\sum_j  Z_j\otimes\gamma_{\fq_\frl}(Z_j)\otimes 1+\sqrt{2}\sum_k T_k\otimes\gamma_{\fq_\frl}(T_k)\otimes 1.
\]
We now recall that
\begin{eqnarray*}
\widehat{D}_{L/L\cap H}(E)&=&-\sum_j  Z_j\otimes\gamma_{\fq_\fl}(Z_j)\otimes 1+\sum_k T_k\otimes\gamma_{\fq_\fl}(T_k)\otimes 1;\\
D_{L/L\cap H}(E)&=&\widehat{D}_{L/L\cap H}(E)-1\otimes\gamma_{\frq_\frl}(c_{\frl,\frl\cap\frh})\otimes 1;\\
D_{L\cap K/L\cap H}(E)&=&-\sum_j  Z_j\otimes\gamma_{\fq_\fl}(Z_j)\otimes 1;\\
D_{L/L\cap K}(\widetilde{E})&=&\sum_k T_k\otimes\gamma_{\fq_\fl}(T_k)\otimes 1.
\end{eqnarray*}
(Note that $D_{L\cap K/L\cap H}(E)=\widehat{D}_{L\cap K/L\cap H}(E)$ since $L\cap K/L\cap H$ is symmetric for triples of type $S$. Also, the definition of this operator has $\gamma_{\frl\cap\frs}$ in place of $\gamma_{\fq_\fl}$, but the two can be identified via the embedding of $C(\frl\cap\frs)$ into $C(\frq_\frl)$.)

It follows that 
\begin{equation}\label{d1}
\imath(D_1)= \sqrt{2}\; \big(D_{L/L\cap H}(E)+1\otimes\gamma_{\frq_\frl}(c_{\frl,\frl\cap\frh})\otimes 1\big)+(1-\sqrt{2})\; D_{L\cap K/L\cap H}(E),
\end{equation}
and also that
\begin{equation}\label{d1 b}
\imath(D_1)= \sqrt{2}\; D_{L/L\cap K}(\widetilde{E})+ D_{L\cap K/L\cap H}(E),
\end{equation}

We next compute the summand $\imath(1\otimes \gamma_\frq(D_2)\otimes 1)$ of $\imath(D_{G/H})$. Since $\rho^-(Z_r)$ and $\rho^-(T_i)$ form an orthonormal basis of $\frq$, with $\|\rho^-(Z_r)\|^2=-1$ and $\|\rho^-(T_i)\|^2=1$, Lemma \ref{alphah2} gives 
\begin{eqnarray*}
\alpha_\fh(\rho^+(T_k))=&-&\sum_{r<s}\underbrace{\lan [\rho^+(T_k),\rho^-(Z_r)],\rho^-(Z_s)\ran }_{=0\text{ since }\fk\perp\fs}\rho^-(Z_r)\rho^-(Z_s)\\
&+& \sum_{r,i}\lan [\rho^+(T_k),\rho^-(Z_r)],\rho^-(T_i)\ran\, \rho^-(Z_r)\rho^-(T_i)\\
&-&\sum_{i<j}\underbrace{\lan [\rho^+(T_k),\rho^-(T_i))],\rho^-(T_j)\ran }_{=0\text{ since }\fk\perp\fs}\rho^-(T_i)\rho^-(T_j)\\
=&&\sum_{r,i}\lan [\rho^+(T_k),\rho^-(Z_r)],\rho^-(T_i)\ran\, \rho^-(Z_r)\rho^-(T_i).
\end{eqnarray*}
It follows that
\eq\label{d2 temp}
D_2=\sum_{k,r,i}\lan [\rho^+(T_k),\rho^-(Z_r)],\rho^-(T_i)\ran\, \rho^-(T_k)\rho^-(Z_r)\rho^-(T_i),
\eeq
and hence
\eq\label{d2 temp 2}
\imath(1\otimes \gamma_\frq(D_2)\otimes 1)=1\otimes \gamma_{\frq_\frl}\Big(\sum_{k,r,i}\lan [\rho^+(T_k),\rho^-(Z_r)],\rho^-(T_i)\ran \, T_k Z_r T_i\Big)\otimes 1.
\eeq
The following lemma will help us to simplify the coefficients in the above expression. 

\begin{lem}\label{omega}
In the above setting, we have
\[
\lan [\rho^+(T_k),\rho^-(Z_r)],\rho^-(T_i)\ran =\lan [T_k,Z_r],T_i\ran .
\]
\end{lem}

We postpone the proof until after we prove Theorem \ref{embeddingthm}. Using Lemma \ref{omega}, we rewrite \eqref{d2 temp 2} as
\eq\label{d2 temp 3}
\imath(1\otimes \gamma_\frq(D_2)\otimes 1)=1\otimes \gamma_{\frq_\frl}\Big(\sum_{k,r,i}\lan [T_k,Z_r],T_i\ran T_k Z_r T_i\Big)\otimes 1.
\eeq
On the other hand, by \eqref{cubicterm} we know
\[
c_{\frl,\frl\cap\frh}=-\Big( \sum_{r<s<t}\lan [Z_r,Z_s],Z_t\ran  Z_r Z_s Z_t + \sum_{r;k<i}\lan [T_k,Z_r],T_i\ran T_k Z_r T_i\Big)
\]
(the terms with $Z_r Z_s T_k$ and $T_i T_j T_k$ vanish because $\frk\perp \frs$). We also know that
\[
c_{\frl\cap\frk,\frl\cap\frh}=- \sum_{r<s<t}\lan [Z_r,Z_s],Z_t\ran  Z_r Z_s Z_t =0
\]
since $L\cap K/L\cap H$ is a symmetric space for triples of type $S$. So
\[
c_{\frl,\frl\cap\frh}=-\sum_{r;k<i}\lan [T_k,Z_r],T_i\ran T_k Z_r T_i,
\]
and comparing this with \eqref{d2 temp 3} we conclude that
\eq\label{d2}
\imath(1\otimes \gamma_\frq(D_2)\otimes 1)=-2\otimes \gamma_{\frq_\frl}(c_{\frl,\frl\cap\frh})\otimes 1.
\eeq

Finally we consider the third summand of $\imath(D_{G/H})$,
\[
\imath(1\otimes D_3)=1\otimes \sum_k\gamma_{\frq_\frl}(T_k)\otimes \beta(\rho^+(T_k)).
\]
We know from Section \ref{triples} that $\rho^+$ maps $\frl\cap\frs$ isometrically onto $\frh\cap\frs$. In particular, $\rho^+(T_k)$ form an orthonormal basis of $\frh\cap\frs$. Also, we can identify
$C(\frq_\frl)=C(\frq_\frl')\bar\otimes C(\frl\cap\frs)$ with $C(\frq_\frl')\bar\otimes C(\frh\cap\frs)$, and $S_{\frq_\frl}=S_{\frq_\frl'}\otimes S_{\frl\cap\frs}$ with $S_{\frq_\frl'}\otimes S_{\frh\cap\frs}$. In this way we see that we can identify
\eq\label{d3}
\imath(1\otimes D_3)=D_{\frh,\frh\cap\frk}(E),
\eeq
with $D_{\frh,\frh\cap\frk}(E)$ acting on $S_{\frh\cap\frs}\otimes E$.

We now add up \eqref{d1}, \eqref{d2} and \eqref{d3} to get the first claim of Theorem \ref{embeddingthm}, or \eqref{d1 b}, \eqref{d2} and \eqref{d3} to get the second claim of Theorem \ref{embeddingthm}.
This finishes the proof of Theorem \ref{embeddingthm} modulo Lemma \ref{omega}.
\epf

\noindent{\it Proof of Lemma \ref{omega}.} Let $\omega$ denote the invariant trilinear alternating form on $\frg$ given by 
\[
\omega(X,Y,Z)=\lan [X,Y],Z\ran .
\]
Suppose $X\in\fl(\nu)$, $Y\in\fl(\nu^\prime)$ and $Z\in\fl(\nu^{\prime\prime})$. We claim that
\eq\label{claim omega}
\omega(\rho^+(X),\rho^-(Y),\rho^-(Z))=\frac{1}{4}d_\nu d_{\nu^\prime} d_{\nu^{\prime\prime}}(1+\nu-\nu^\prime-\nu^{\prime\prime})\ \omega(X,Y,Z).
\eeq
If we prove this claim, the lemma follows immediately by applying the claim to $X=T_k$, $Y=Z_r$, $Z=T_i$, since in that case 
\[
\nu=\nu''=0;\quad \nu'=-1;\qquad d_\nu=d_{\nu''}=\sqrt{2};\quad d_{\nu'}=1.
\]
To prove \eqref{claim omega}, we first recall that 
\[
\rho^+(X)=\frac{d_\nu}{2}(X+\sigma(X)),\quad\rho^-(Y)=\frac{d_{\nu'}}{2}(Y-\sigma(Y)),\quad\rho^-(Z)=\frac{d_{\nu''}}{2}(Z-\sigma(Z)),
\]
so \eqref{claim omega} will follow if we prove
\eq\label{claim omega2}
\omega(X+\sigma(X),Y-\sigma(Y),Z-\sigma(Z))=2\omega(X,Y,Z).
\eeq
By definition of $\frl(\nu)$,
\begin{multline}\label{om3}
\omega(\sigma(X),Y,Z)=\lan[\sigma(X),Y],Z\ran=\\ 
\lan\sigma(X),\underbrace{[Y,Z]}_{\in\fl}\ran=
\nu\lan X,[Y,Z]\ran=\nu\,\omega(X,Y,Z).
\end{multline}
Since $\omega$ is skew symmetric, this implies that also
\eq\label{om4}
\omega(X,\sigma(Y),Z)=\nu'\,\omega(X,Y,Z),\qquad \omega(X,Y,\sigma(Z))=\nu''\,\omega(X,Y,Z).
\eeq
The $\sigma$-invariance of $\omega$ now implies
\eq\label{om5}
\omega(\sigma(X),\sigma(Y),Z)=\omega(X,Y,\sigma(Z))=\nu^{\prime\prime}\,\omega(X,Y,Z),
\eeq
so by the skew symmetry of $\omega$ also 
\eq\label{om6}
\omega(\sigma(X),Y,\sigma(Z))=\nu^{\prime}\,\omega(X,Y,Z),\qquad  \omega(X,\sigma(Y),\sigma(Z))=\nu\,\omega(X,Y,Z).
\eeq
The equation \eqref{claim omega2} follows immediately from \eqref{om3} -- \eqref{om6}.
\qed

\begin{rem} Similar formulas for $\imath(D_{G/H}(E))$ can also be obtained for triples $(G,H,L)$ of type $I$ or $N$.
They are however more complicated, and therefore less likely to have nice applications. Especially, we do not have a nice interpretation as above for the operator $D_3$ and its action seems hard to understand.

If however the subgroup $H$ is simply connected, the spin module $S_\fq$ admits an action of $H$ and one can choose the $\fh$-module $E$ to be the trivial representation. In this case, the operator $D_3$ is zero and the above problem with interpretation disappears. 
\end{rem}

\section{Example: $(G,H,L)=(G'\times G',\Delta(G'\times G'),G'\times K')$ for\\ $G'=SL(2,\bbR)$ and $K'=SO(2)$}\label{sec sl2}

Let $(G,H,L)$ be a triple $(G'\times G',\Delta(G'\times G'),G'\times K')$, where $G'$ is a non-compact connected semisimple real Lie group with finite center, with maximal compact subgroup $K'$, $G=G'\times G'$, $H=\Delta(G'\times G')$ and $L=G'\times K'$.
Note that the maximal compact subgroup of $G$ is $K=K'\times K'$. Moreover, $L\cap H=\Delta(K'\times K')$, and $L\cap K=K$.

If $\frg'=\frk'\oplus\frs'$ is the Cartan decomposition of the Lie algebra of $G'$, then we have
\[
\frg=\frg'\times \frg',\quad \frh=\Delta(\frg'\times \frg'),\quad \frl=\frg'\times \frk',
\]
and furthermore
\begin{gather*}
\frq=\Delta^-(\frg'\times \frg'),\quad \frq_\frl=\frq_\frl'\oplus(\frl\cap\frs)=\Delta^-(\frk'\times\frk')\oplus(\frs'\times 0);\\
\frh\cap\frk=\Delta(\frk'\times\frk')=\frl\cap\frh,\quad \frh\cap\frs=\Delta(\frs'\times\frs');\\ \frq\cap\frk=\Delta^-(\frk'\times\frk'),\quad\frq\cap\frs=\Delta^-(\frs'\times\frs');\\
\frl\cap\frk=\frk,\quad\frl\cap\frs=\frs'\times 0
\end{gather*}
where $\Delta^-$ denotes the anti-diagonal. From now on we let $G'=SL(2,\bbR)$ and $K'=SO(2)$. Let
\[
h=\begin{pmatrix} 0&-i\cr i&0\end{pmatrix},\qquad 
e=\frac{1}{2}\begin{pmatrix} 1&i\cr i&-1\end{pmatrix},\qquad
f=\frac{1}{2}\begin{pmatrix} 1&-i\cr -i&-1\end{pmatrix}
\]
be the usual $\frsl(2)$-basis of $\frg'_\bbC=\frsl(2,\bbC)$, with commutators
\[
[h,e]=2e,\qquad [h,f]=-2f,\qquad [e,f]=h.
\]
We consider the trace form $\lara$ on $\frg'$ and the direct sum form on $\frg$. It is easy to check that
\[
\htil=\frac{i}{\sqrt{2}}h,\quad \etil=\frac{1}{\sqrt{2}}(e+f),\quad \ftil=\frac{i}{\sqrt{2}}(e-f)
\]
form an orthonormal basis of $\frg'$, i.e., they are orthogonal to each other and satisfy
\[
\lan \htil,\htil\ran=-1,\quad \lan \etil,\etil\ran=\lan \ftil,\ftil\ran =1.
\]
So if we set 
\[
Z=\frac{1}{\sqrt{2}}(\htil,-\htil),\quad T_1=(\etil,0),\quad T_2=(\ftil,0),
\]
we see that $Z$ forms an orthonormal basis of $\frq_\frl'$ (i.e., $\lan Z,Z\ran=-1$), while $T_1$ and $T_2$ form an orthonormal basis of $\frl\cap\frs$
(i.e., $\lan T_1,T_1\ran=\lan T_2,T_2\ran=1$, $\lan T_1,T_2\ran=0$).

Since $L\cap K$ is $\sigma$-stable, then $\la=-1$, i.e., $Z\in\frq$ and $d_\lambda=1$. Moreover, $\fl\cap\fs$ is irreducible under the action of $L\cap H$, with $\mu=0$ and $d_\mu=\sqrt{2}$. We see immediately that
\begin{gather}
\label{rho pm}
\rho^+(Z)=0,\qquad\rho^-(Z)=Z;\\
\label{rho pm2}
\rho^+(T_1)=\frac{1}{\sqrt{2}}(\etil,\etil),\qquad\rho^-(T_1)=\frac{1}{\sqrt{2}}(\etil,-\etil);\\
\label{rho pm3}
\rho^+(T_2)=\frac{1}{\sqrt{2}}(\ftil,\ftil),\qquad\rho^-(T_2)=\frac{1}{\sqrt{2}}(\ftil,-\ftil).
\end{gather}
Since $SL(2,{\mathbb R})$ acts on any finite-dimensional representation of ${\mathfrak s}{\mathfrak l}(2,{\mathbb C})$, we may choose $E$ to be any finite-dimensional representation of $\Delta({\mathfrak s}{\mathfrak l}(2,{\mathbb C})\times {\mathfrak s}{\mathfrak l}(2,{\mathbb C}))$. By the second statement of Theorem \ref{embeddingthm}, the embedding of $D_{G/H}(E)$ can be written as
\eq
\label{emb D}
\imath(D_{G/H}(E))=\sqrt{2}D_{L/L\cap K}(\widetilde{E})+D_{L\cap K/L\cap H}(E)-2\otimes\gamma_{\frq_\frl}(c_{\frl,\frl\cap\frh})\otimes 1+1\otimes D_{\frh,\frh\cap\frk}(E).
\eeq
This is an equality of operators acting on the space
\eq
\label{sections}
C^\infty(L,S_{\frl\cap\frs}\otimes C^\infty(L\cap K,S_{\frq_\frl'}\otimes E)^{L\cap H})^{L\cap K}.
\eeq
The $L\cap H$-invariants are taken with respect to the right translation on $C^\infty(L\cap K)$ and lift of $\fl\cap\fh$-action on $S_{\fq_\fl^\prime}\otimes E$ 
and the $L\cap K$-invariants are taken with respect to the right translation on $C^\infty(L)$ and lift of $\fl\cap\fk$-action on $S_{\fl\cap\fs}\otimes \widetilde{E}$. 
Recall that we are using the isometry $\rho^+:\frl\cap\frs\to\frh\cap\frs$ to identify $S_{\frl\cap\frs}$ with $S_{\frh\cap\frs}$.

We first consider the action of the summand
\[
-2\otimes\gamma_{\frq_\frl}(c_{\frl,\frl\cap\frh})\otimes 1
\]
of \eqref{emb D}. By \eqref{cubicterm}, 
\[
c_{\frl,\frl\cap\frh}=-ZT_1T_2,
\]
so we have to determine the action of $2ZT_1T_2$ on the spin module $S_{\frq_\frl}=S_{\frl\cap\frs}\otimes S_{\frq_\frl'}$. To construct the spin module we use the dual isotropic vectors
\[
u=\frac{1}{\sqrt{2}}(T_1+iT_2),\qquad v=\frac{1}{\sqrt{2}}(T_1-iT_2).
\]
Then $S_{\frq_\frl}$ is spanned by 1 and $u$, and the action of $u$ and $v$ is given by
\[
u\cdot 1=u,\quad u\cdot u=0;\qquad v\cdot 1=0,\quad v\cdot u=1.
\]
Since $T_1=\frac{1}{\sqrt{2}}(u+v)$ and $T_2=-\frac{i}{\sqrt{2}}(u-v)$, they act by 
\eq
\label{spin ls action}
T_1\cdot 1=\frac{1}{\sqrt{2}}u,\quad T_1\cdot u=\frac{1}{\sqrt{2}};\qquad T_2\cdot 1=-\frac{i}{\sqrt{2}}u,\quad T_2\cdot u=\frac{i}{\sqrt{2}}.
\eeq
Finally, $Z$ acts by $\pm i/\sqrt{2}$ on $S^\pm$, so
\[
Z\cdot 1= \frac{i}{\sqrt{2}},\quad Z\cdot u =-\frac{i}{\sqrt{2}}u.
\]
A straightforward computation now proves
\begin{lem}
\label{cubic action sl2}
The summand 
\[
-2\otimes\gamma_{\frq_\frl}(c_{\frl,\frl\cap\frh})\otimes 1 =1\otimes \gamma_{\frq_\frl}(2ZT_1T_2)\otimes 1
\]
of the expression for $\imath(D_{G/H}(E))$ given by \eqref{emb D} acts by the scalar $\frac{1}{\sqrt{2}}$ on the whole space \eqref{sections}.
\end{lem}

To understand the action of the other summands of \eqref{emb D}, we describe
\[
\widetilde E = C^\infty(L\cap K,S_{\frq_\frl'}\otimes E)^{L\cap H}
\]
more explicitly. The group 
\[
L\cap K=K'\times K'=SO(2)\times SO(2)
\] 
is abelian, compact and connected, and its Lie algebra is
\[
\frl\cap\frk=\frl\cap\frh\oplus \frq_\frl'.
\]
The Lie algebra $\frl\cap\frh$ is one-dimensional, spanned by
\[
W=\frac{1}{\sqrt{2}}(\htil,\htil).
\]
Furthermore, $W$ and $Z$ form an orthonormal basis for the abelian Lie algebra $\frl\cap\frk$. 
We denote by $\bbC_{a,b}$ the one-dimensional character of $L\cap K$ on which $(h,0)$ acts by the scalar $a$ and $(0,h)$ acts by the scalar $b$. Since $h\in\frk'_\bbC$ acts by integers on the characters of $K'=SO(2)$, $a$ and $b$ are integers. Furthermore, $W=\frac{i}{2}(h,h)$ acts on 
$\bbC_{a,b}$ by the scalar $\frac{i(a+b)}{2}$, while $Z=\frac{i}{2}(h,-h)$ acts on $\bbC_{a,b}$ by the scalar $\frac{i(a-b)}{2}$.

Since the space of smooth vectors in $L^2(L\cap K)$ coincides with $C^\infty(L\cap K)$ and the dual of $\bbC_{a,b}$ is $\bbC_{-a,-b}$, Peter-Weyl theorem implies that
\[
C^\infty(L\cap K)=\overline{\bigoplus}_{a,b}\bbC_{-a,-b}\otimes\bbC_{a,b}
\]
with the left regular action of $L\cap K$ realized on the first factor of each summand, and the right regular action of $L\cap K$ realized on the second factor of each summand. 
Here $\overline{\bigoplus}$ denotes the closure of the algebraic direct sum in the smooth topology (see \cite{MO21} for more detail in a more general setting.) Since the $L\cap H$-invariants are taken with respect to the right regular action, we have (up to passing to a dense subspace)
\eq
\label{dec sec}
C^\infty(L\cap K,S_{\frq_\frl'}\otimes E)^{L\cap H}=\bigoplus_{a,b}\bbC_{-a,-b}\otimes(\bbC_{a,b}\otimes S_{\frq_\frl'}\otimes E)^{L\cap H}.
\eeq
Since $L\cap H$ is connected, and $\frl\cap\frh$ is spanned by $(h,h)$, the $L\cap H$-invariants are simply the 0-eigenspace of $(h,h)$. The map
\[
\alpha_{\frl\cap\frk,\frl\cap\frh}:\frl\cap\frh\to C(\frq_\frl')
\]
is zero since $\frl\cap\frk$ is abelian, so $(h,h)$ acts by $0$ on $S_{\frq_\frl'}$. Furthermore, $(h,h)$ acts by $a+b$ on $\bbC_{a,b}$. It follows that taking the $L\cap H$-invariants in \eqref{dec sec} amounts to pairing each $\bbC_{a,b}$ with the $(-a-b)$-weight space of $E$. So \eqref{dec sec} can be
rewritten as 
\eq
\label{dec sec 2}
C^\infty(L\cap K,S_{\frq_\frl'}\otimes E)^{L\cap H}=\bigoplus_{a,b;\,a+b\in\Delta(E)}\bbC_{-a,-b}\otimes\bbC_{a,b}\otimes S_{\frq_\frl'}\otimes E_{-a-b},
\eeq
with $\Delta(E)$ denoting the set of weights of $E$.
So the space \eqref{sections} can be identified with
\eq
\label{sections 2}
C^\infty\big[L,S_{\frl\cap\frs}\otimes \big(\bigoplus_{a,b;\,a+b\in\Delta(E)}\bbC_{-a,-b}\otimes\bbC_{a,b}\otimes S_{\frq_\frl'}\otimes E_{-a-b}\big)\big]^{L\cap K}.
\eeq
We now consider the summand 
\eq\label{d lk lh}
D_{L\cap K/ L\cap H}(E)=-Z\otimes Z
\eeq 
of \eqref{emb D}. Recall that the first factor $Z$ is acting 
on $C^\infty(L\cap K,S_{\frq_\frl'}\otimes E)^{L\cap H}$
by the right regular action, while the second factor $Z$ is acting on $S_{\frq_\frl'}$. So the action of $D_{L\cap K/ L\cap H}(E)$ on the space \eqref{sections 2} is 
on $\bbC_{a,b}\otimes S_{\frq_\frl'}$.

Since the spin module $S_{\frq_\frl'}$ is one-dimensional, with $Z\in C(\frq_\frl')$ acting by the scalar $\frac{i}{\sqrt{2}}$, and since $Z$ operates on $\bbC_{a,b}$ by the scalar $\frac{i(a-b)}{2}$, 
the operator $D_{L\cap K/ L\cap H}(E)$ given by \eqref{d lk lh} acts on $\bbC_{a,b}\otimes S_{\frq_\frl'}$ by the scalar $\frac{a-b}{2\sqrt{2}}$. In this way we have diagonalized the action of $D_{L\cap K/ L\cap H}(E)$ on the space \eqref{sections 2}: the eigenvalues are 
\[
\frac{k}{2\sqrt{2}},\qquad k\in\bbZ,
\] 
and the eigenspace corresponding to $k$ is the space
\eq\label{e spaces}
C^\infty\big[L,S_{\frl\cap\frs}\otimes \big(\bigoplus_{a,b;\,a-b=k,\,a+b\in\Delta(E)}\bbC_{-a,-b}\otimes\bbC_{a,b}\otimes S_{\frq_\frl'}\otimes E_{-a-b}\big)\big]^{L\cap K}.
\eeq
In particular, if we set $k=-2$ and recall Lemma \ref{cubic action sl2}, we obtain
\begin{cor}
\label{cubic plus lk lh}
The part 
\[
D_{L\cap K/L\cap H}(E)-2\otimes\gamma_{\frq_\frl}(c_{\frl,\frl\cap\frh})\otimes 1
\]
of the expression for $\imath(D_{G/H}(E))$ given by \eqref{emb D} acts by 0 on the subspace
\[
C^\infty\big[L,S_{\frl\cap\frs}\otimes \big(\bigoplus_{a,b;\,a-b=-2,\,a+b\in\Delta(E)}\bbC_{-a,-b}\otimes\bbC_{a,b}\otimes S_{\frq_\frl'}\otimes E_{-a-b}\big)\big]^{L\cap K}.
\]
of the space \eqref{sections 2}. This subspace is nonzero if and only if $E$ has even weights.
\end{cor}

We now turn our attention to the summand $D_{\frh,\frh\cap\frk}(E)$ of the expression for $\imath(D_{G/H}(E))$ given by \eqref{emb D}. Identifying $\frh$ with $\frg'=\frsl(2,\bbR)$, we can use the basis $e,f,h$ of $\frg'_\bbC$ to write 
\[
D_{\frh,\frh\cap\frk}(E)=e\otimes f+f\otimes e
\]
with the second factor acting on 
\[
S_{\frh\cap\frs}= \bbC\,1\oplus\bbC\,e
\]
in the usual way. If $E$ has highest weight $2m$, $m\in\bbZ_+$, then it is well known and easy to see that
\[
\ker D_{\frh,\frh\cap\frk}(E) = (E_{2m}\otimes e)\oplus (E_{-2m}\otimes 1).
\]
Identifying $S_{\frh\cap\frs}$ with $S_{\frl\cap\frs}$, the basis of the spin module becomes $\{1,u\}$, and we conclude
\begin{cor}
\label{D gh D l lk}
If the highest weight of $E$ is $2m$, $m\in\bbZ_+$, 
then the operators $\imath(D_{G/H}(E))$ and $\sqrt{2}D_{L/L\cap K}(\widetilde{E})$ are equal on the subspace
\begin{multline*}
C^\infty\big(L,u\otimes \bbC_{m+1,m-1}\otimes\bbC_{-m-1,-m+1}\otimes S_{\frq_\frl'}\otimes E_{2m}\,{\textstyle{\bigoplus}}\\ 1\otimes \bbC_{-m+1,-m-1}\otimes\bbC_{m-1,m+1}\otimes S_{\frq_\frl'}\otimes E_{-2m}\big)^{L\cap K}
\end{multline*}
of the space \eqref{sections 2}. In particular, $\imath(D_{G/H}(E))$ and $D_{L/L\cap K}(\widetilde{E})$ have the same kernel on this subspace.
\end{cor}

We now discuss the kernel of $D_{L/L\cap K}(\widetilde{E})$ on the subspace of the space \eqref{sections 2} described in Corollary \ref{D gh D l lk}. We first recall that in the definition of the space 
\[
C^\infty(L,S_{\frl\cap\frs}\otimes\widetilde E)^{L\cap K}\cong (C^\infty(L)\otimes S_{\frl\cap\frs}\otimes\widetilde E)^{L\cap K},
\] 
the $L\cap K$-invariants are taken with respect to the right regular action on $C^\infty(L)$, the spin action on $S_{\frl\cap\frs}$, and the left regular action on
\[
\widetilde E=C^\infty(L\cap K,S_{\frq_\frl'}\otimes E)^{L\cap H}.
\]
This means that we can identify the subspace described in Corollary \ref{D gh D l lk} with
\eq
\label{sec 3}
C^\infty\big(L,u\otimes \bbC_{m+1,m-1}\,{\textstyle{\bigoplus}}\, 1\otimes \bbC_{-m+1,-m-1}\big)^{L\cap K};
\eeq 
the omitted factors are one-dimensional and carry no relevant action. 

We now use $L=G'\times K'$ and $L\cap K=K'\times K'$ to rewrite the space \eqref{sec 3} as
\begin{multline}
\label{sec 5}
C^\infty\big(G',u\otimes \bbC_{m+1})^{K'}\otimes C^\infty\big(K',\bbC_{m-1})^{K'}
{\textstyle{\bigoplus}}\\ 
C^\infty(G',1\otimes\bbC_{-m+1})^{K'}\otimes C^\infty\big(K',\bbC_{-m-1})^{K'}.
\end{multline} 
Here and below for any $k\in\bbZ$, $\bbC_k$ denotes the character of $K'$ on which $h$ acts by $k$.
Peter-Weyl Theorem for $K'$ implies that for any such character $\bbC_k$, 
\[
C^\infty\big(K',\bbC_k)^{K'}=\overline{\bigoplus}_{n\in\bbZ}\bbC_{-n}\otimes(\bbC_n\otimes\bbC_k)^{K'}=\bbC_k.
\]
So \eqref{sec 5} becomes
\eq
\label{sec 6}
C^\infty\big(G',u\otimes\bbC_{m+1})^{K'}\otimes \bbC_{m-1}\,{\textstyle{\bigoplus}}\, C^\infty(G',1\otimes\bbC_{-m+1})^{K'}\otimes\bbC_{-m-1}.
\eeq
Note that $D_{L/L\cap K}(\widetilde{E})=D_{G'/K'}$ is acting on the first factors in \eqref{sec 6}. To determine the kernel of this operator, we use a reciprocity result of \cite{MZ14}. More precisely, let $V$ be a smooth admissible representation of $G'$, $V_{K'}$ the space of $K'$-finite vectors and $V_{K'}^*$ the $K'$-finite dual of $V_{K'}$. Let $S_{\fs'}$ be the spin module for $\fk'$ and $F$ a finite dimensional representation of $\fk'$ such that $S_{\fs'}\otimes F$ lifts to a $K'$-representation. Proposition 1.8 of \cite{MZ14} says that the map
$$\Psi:\text{Hom}_{G'}(V,C^\infty(G'/K',{\mathcal S_{\fs'}}\otimes{\mathcal F}))\rightarrow \text{Hom}_{K'}(F^*,S_{\fs'}\otimes V_{K'}^*)$$
defined by
\[
\Psi(T)(f^*)(v)=(\id_{S_{\frs'}}\otimes f^*) [T(v)(1)],
\]
where we identify $S_{\fs'}\otimes V_{K'}^*$ with $\Hom(V_{K'},S_{\frs'})$, induces an isomorphism
\eq
\label{reciprocity}
\Hom_{G'}(V,\text{ker}(D_{G'/K'}(F)))\simeq\Hom_{K'}(F^*,\text{ker}(D_{V_{K'}^*}))
\eeq
where $F^*$ is the dual of $F$. 
This isomorphism does not involve the spin module, so we can replace the spin module by any of its irreducible components and still have the same statement. This means that the kernel of $D_{G'/K'}$ on $C^\infty\big(G',u\otimes\bbC_{m+1})^{K'}$ contains irreducible representations $(\pi,V_\pi)$ of $G'$ with the property that the kernel of $D_{\frg',\frk'}$ on $u\otimes V_\pi^*$ contains $\bbC_{m+1}^*=\bbC_{-m-1}$. There is exactly one such $\pi^*$: the irreducible representation with 
highest weight $-m-2$. So $\pi$ is the unique irreducible representation with lowest weight $m+2$, i.e., the discrete series representation with lowest weight $m+2$ (recall that $m\geq 0$, so $m+2\geq 2$). We denote this representation by $DS^+_{m+2}$.

Analogously, the kernel of $D_{G'/K'}$ on $C^\infty\big(G',1\otimes\bbC_{-m+1})^{K'}$ contains irreducible representations $(\pi,V_\pi)$ of $G'$ with the property that the kernel of $D_{\frg',\frk'}$ on $1\otimes V_\pi^*$ contains $\bbC_{-m+1}^*=\bbC_{m-1}$. There is exactly one such $\pi^*$: the irreducible representation with 
lowest weight $m$. The corresponding $\pi$ is the irreducible representation with highest weight $-m$, so it is
the discrete series representation $DS^-_{-m}$ with highest weight $-m$ if $m\geq 2$, the limit of discrete series representation $LDS^-$ with highest weight $-1$ if $m=1$, or the trivial representation $\bbC$ if $m=0$.

The above discussion and \eqref{sec 6} imply 

\begin{thm}
\label{thm sl2}
Suppose that the highest weight of $E$ is $2m$, $m\in\bbZ_+$. Then the kernel of $D_{L/L\cap K}(\widetilde{E})$ on $C^\infty(L,S_{\frl\cap\frs}\otimes\widetilde E)^{L\cap K}$ contains the following representations of $L=G'\times K'$:
\begin{gather*}
(DS_{m+2}^+\otimes \bbC_{m-1})\oplus (DS_{-m}^-\otimes \bbC_{-m-1})\qquad\text{if}\quad m\geq 2; \\
(DS_3^+\otimes \bbC_0)\oplus (LDS^-\otimes \bbC_{-2})\qquad\text{if}\quad m=1; \\
(DS^+_{2}\otimes\bbC_{-1})\oplus (\bbC\otimes\bbC_{-1})\qquad\text{if}\quad m=0.
\end{gather*}
\end{thm}

It follows that the same $L$-representations appear in the kernel of $\imath(D_{G/H})$, and thus also in the kernel of $D_{G/H}$. On the other hand, the kernel of $D_{G/H}$ is a representation of $G$ and not only of $L$, so our $L$-representations inside $\ker(D_{G/H})$ generate $G$-representations inside $\ker(D_{G/H})$. t is possible (but quite complicated) to identify the representations of $G$ in $\ker(D_{G/H})$ whose restrictions to $L$ contain the representations described in Theorem \ref{thm sl2}.

\end{document}